\documentclass[11pt]{amsart}
\usepackage{amsmath,amsfonts,latexsym,amssymb,amsthm,mathrsfs,array}

\usepackage{amsfonts}
\usepackage[usenames]{color}
\usepackage[colorlinks=true,linkcolor=webgreen,filecolor=webbrown,citecolor=webgreen]
{hyperref}
\definecolor{webgreen}{rgb}{0,.5,0}
\definecolor{webbrown}{rgb}{.6,0,0}

\def\hh{{\mathfrak{h}}}

\newtheorem{dfn}{Definition}[section]
\newcommand{\bdfn}{\begin{dfn}\rm}
\newcommand{\edfn}{\end{dfn}}
\newtheorem{thm}[dfn]{Theorem}
\newcommand{\bthm}{\begin{thm}}
\newcommand{\ethm}{\end{thm}}
\newtheorem{lmma}[dfn]{Lemma}                   
\newcommand{\blmma}{\begin{lmma}}                   
\newcommand{\elmma}{\end{lmma}}                   
\newtheorem{ppsn}[dfn]{Proposition}
\newcommand{\bppsn}{\begin{ppsn}}
\newcommand{\eppsn}{\end{ppsn}}
\newtheorem{crlre}[dfn]{Corollary}
\newcommand{\bcrlre}{\begin{crlre}} 
\newcommand{\ecrlre}{\end{crlre}}
\newtheorem{rmk}[dfn]{Remark}
\newcommand{\brmk}{\begin{rmk}\rm} 
\newcommand{\ermk}{\end{rmk}}
\newtheorem{conj}[dfn]{Conjecture}
\newcommand{\bconj}{\begin{conj}\co} 
\newcommand{\econj}{\end{conj}}

\numberwithin{equation}{section}

\title{Generalized Casimir Operators for Lie Superalgebras}
\author{S. Eswara Rao}
\address{S. Eswara Rao, School of Mathematics, Tata Institute of Fundamental Research, Homi Bhabha Road, Colaba, Mumbai 400005, India.}                  
\email{senapati@math.tifr.res.in, sena98672@gmail.com}
\date{}

\begin{document} 
\maketitle
\begin{abstract}
In this paper, we define generalized Casimir operators for a loop contragredient Lie superalgebra and prove that they commute with the underlying Lie superalgebra. These operators have applications in the decomposition of tensor product modules. We further introduce the notion of generalized Gelfand invariants for the loop general linear Lie superalgebra and show that they also commute with the underlying Lie superalgebra. These operators when applied to a highest weight vector in a tensor product module again induces a new highest weight vector.\\\\           
{\bf{MSC}(2010):} 17B65, 17B70.\\
{\bf{KEY WORDS}:} Contragredient Lie Superalgebras, Generalized Casimir Operators, Generalized Gelfand Invariants, Evaluation Modules, Tensor Product Modules.      
\end{abstract}

\section{Introduction}
The central purpose of this paper is to construct generalized Casimir operators for loop contragredient Lie superalgebras. A contragredient Lie superalgebra $\mathfrak{g}$ is defined by means of a complex matrix (see Section \ref{S2}). For any associative, commutative and finitely generated algebra $B$ with unity, we first consider the loop contragredient Lie superalgebra $\mathfrak{g} \otimes B$. We then define generalized Casimir operators in the completion of its universal enveloping algebra $\overline{U}(\mathfrak{g} \otimes B)$. Henceforth we prove that these operators commute with the action of $\mathfrak{g} \cong \mathfrak{g} \otimes \mathbb{C}$ (see Theorem \ref{T3.2}). These generalized Casimir operators play an important role in the tensor product decomposition of $\mathfrak{g}$-modules (see Section \ref{S4}).

In recent times, loop algebras of basic classical Lie superalgebras have been studied. The finite dimensional irreducible modules for the loop algebra of basic classical Lie superalgebras are classified in \cite{RZ} and \cite{R1} when $B$ is a Laurent polynomial algebra in finitely many variables. The twisted case was settled by A. Savage in \cite{S}.

In the current paper, we work with a symmetrizable contragredient Lie superalgebra which admits a non-degenerate, super symmetric, invariant and even bilinear form (see (2.8)). Using this bilinear form, we introduce Casimir operators for $\mathfrak{g}$. In the literature, there are two kinds of Casimir operators-one is defined by Wakimoto \cite[Theorem 2.17]{W} (page-97) and the other is given by Musson \cite[Section 18.4]{M} (page-414). Our operators coincide with the ones mentioned in \cite{W}, but here we deal with much more general Lie superalgebras than in \cite{W}.

We now describe the contents of our paper. In Section \ref{S3}, we recall the basic definition of symmetrizable contragredient Lie superalgebra. They admit non-degenerate bilinear forms satisfying (2.8). Using this form, we define Casimir operators for $\mathfrak{g}$. The proof that they commute with $\mathfrak{g}$ is postponed to Section \ref{S3} (see Theorem \ref{T3.2}) where a much more general result is proved. In the remaining portion of this section, we work out Casimir operators for rank one Lie superalgebras and verify that they are all central.

In Section \ref{S3}, we define the loop contragredient Lie superalgebra $\mathfrak{g} \otimes B$ and the generalized Casimir operators. In Theorem \ref{T3.2}, we establish that these operators commute with $\mathfrak{g} \cong \mathfrak{g} \otimes \mathbb{C}$. In the rest of this section, we work out generalized Casimir operators for rank one Lie superalgebras and check that they commute with the action of $\mathfrak{g}$. These verifications are utilized in the proof of Theorem \ref{T3.2}.

In the last section, we specialize to the case where $\mathfrak{g}= \mathfrak{gl}(M,N)$. We first define the Gelfand invariants which are central
and are known to generate the center of $U(\mathfrak{g})$ as an algebra. Similar to the previous section, for each Gelfand invariant $T_k$($k \geqslant 1$) and for $a_1, \cdots, a_k \in B$, we define $T_k(a_1, \cdots, a_k)$ in the universal enveloping algebra $U(\mathfrak{g} \otimes B)$ and verify that it commutes with the action of $\mathfrak{g}$ (see Lemma \ref{L4.2}).   

The rest of Section \ref{S4} primarily focuses on evaluation modules for $\mathfrak{g} \otimes B$ with $B= \mathbb{C}[t,t^{-1}]$ which is the Laurent polynomial algebra in one variable (see \cite{Ku,R2} for the Lie algebra case). We closely follow these references. Evaluation modules are actually modules for finitely many copies of $\mathfrak{g}$ and thereby lifted to $\mathfrak{g} \otimes B$ by an evaluation map. In this case, there exists a co-finite ideal $J$ of $B$ such that $\mathfrak{g} \otimes J$ acts trivially on the module. Thus for each $k \geqslant 1$, there are only finitely many operators $T_k(a_1, \cdots, a_k)$ where $a_1, \cdots, a_k$ runs over a basis of $B/J$. These operators are very effective for understanding tensor product decompositions of $\mathfrak{g}$-modules. If we take a highest weight vector in the tensor product and apply these operators, we usually obtain a new highest weight vector. In the setup of Lie algebras, the action of these operators is irreducible (see \cite{Ku,R2}). In \cite{G}, Gorelik defined anti-invariants for Lie superalgebras. In case of general linear superalgebras, we provide a class of anti-invariants (see Section \ref{S5}).  
       
\section{Contragredient Lie Superalgebras and Casimir Operator}\label{S2}  
Throughout this paper, all the vector spaces, algebras and tensor products are over the field of complex numbers $\mathbb{C}$. We shall denote the set of all non-negative integers by $\mathbb{Z}_{+}$, the set of natural numbers by $\mathbb{N}$ and the cyclic group of two elements $\{\overline{0}, \overline{1} \}$ by $\mathbb{Z}_2$.   
   
A Lie superalgebra is a $\mathbb{Z}_2$-graded vector space $\mathfrak{g}= \mathfrak{g}_{\overline{0}} \bigoplus \mathfrak{g}_{\overline{1}}$ equipped with a $\mathbb{C}$-bilinear form $[\cdot, \cdot] : \mathfrak{g} \times \mathfrak{g} \longrightarrow \mathfrak{g}$ called the Lie superbracket, satisfying the following conditions--      
\begin{enumerate}
\item $[\mathfrak{g}_{\overline{i}}, \mathfrak{g}_{\overline{j}}] \subseteq  
\mathfrak{g}_{\overline{i} + \overline{j}}$ , 
\item $[x,y]=-(-1)^{\overline{i} \overline{j}} [y,x]$ , 
\item $\big [[x,y],z \big ] = \big [x, [y,z] \big ] - (-1)^{\overline{i} \overline{j}} \big [y, [x,z] \big ]$                           	
\end{enumerate}                     
for all homogeneous elements $x \in \mathfrak{g}_{\overline{i}}, y \in \mathfrak{g}_{\overline{j}}$ and $z \in \mathfrak{g}$.\\  
We shall also denote by $|x|$ the degree of the homogeneous element $x$.\\
In this section we recall the definition of contragredient Lie superalgebra from Musson's book \cite{M}. The main aim of this section is to define the Casimir operator and establish that it always commutes with the universal enveloping algebra (UEA) of this contragredient Lie superalgebra. The contragredient Lie superalgebras are much more general than BKM algebras as defined in Wakimoto's book \cite{W}. BKM algebras do not cover most of the simple finite dimensional Lie superalgebras. We work with symmetrizable Cartan matrices so that the corresponding contragredient Lie superalgebras admit non-degenerate super symmetric invariant bilinear form. Using this bilinear form, we define Casimir operator in the completion of UEA. The operator is similar to the one defined by Wakimoto \cite[Theorem 2.17]{W} (page-97), but defers from the one in Musson's book \cite[Section 18.4]{M} (page-414).  
                 
(2.1) To define contragredient Lie superalgebra, we start with a finite set $\text{I}= \{1, \cdots, l \}$ and a subset $\tau \subseteq \text{I}$, a complex matrix $A=(a_{ij})_{i,j \in \text{I}}$ and a complex vector space $\mathfrak{h}$ of dimension $l +$ co-rank$A$ together with its dual space $\mathfrak{h}^{*}$. There exists linearly independent sets $\pi= \{\alpha_1, \cdots, \alpha_l \} \subseteq \mathfrak{h}^{*}$ and $\pi^{\vee}= \{\alpha_{1}^{\vee}, \cdots, \alpha_{l}^{\vee} \} \subseteq \mathfrak{h}$ such that $\alpha_j(h_i)=a_{ij}$. Let $\widetilde{\mathfrak{g}}(A, \tau)$ be a Lie superalgebra generated by $\mathfrak{h}$ and $\{e_i,f_i \}_{i \in \text{I}}$ with defining relations given by       
\begin{enumerate}
\item $[e_i,f_j]= \delta_{ij}\alpha_{i}^{\vee}$,
\item $[h,e_j]= \alpha_j(h)e_j$,
\item $[h,f_j]= -\alpha_j(h)f_j$,
\item $[h,h^{\prime}]=0 \ \forall \ h,h^{\prime} \in \mathfrak{h}$.        
\end{enumerate}
We define degree($e_i$) = degree($f_i$) = $\overline{1}$, if $i \in \tau$ along with \\ 
degree($e_i$) = degree($f_i$) = $\overline{0}$, if $i \notin \tau$
and degree($h$) = $0 \ \forall \ h \in \mathfrak{h}$.\\ 
The above grading generates a $\mathbb{Z}_2$-grading on       
\begin{align*}
\widetilde{\mathfrak{g}}(A, \tau) = \widetilde{\mathfrak{g}}(A, \tau)_{\overline{0}} \oplus \widetilde{\mathfrak{g}}(A, \tau)_{\overline{1}}.         
\end{align*}

(2.2) Let $\widetilde{N}^{+}$ and $\widetilde{N}^{-}$ be the Lie subalgebras generated by $e_1, \cdots e_l$ and $f_1, \cdots, f_l$ respectively. Following \cite{M}, let us decompose   
\begin{align*}
\widetilde{\mathfrak{g}}(A, \tau) =  \widetilde{N}^{-} \oplus \mathfrak{h} \oplus \widetilde{N}^{+}.            
\end{align*}
Among all the ideals of $\widetilde{\mathfrak{g}}(A, \tau)$ intersecting 
$\mathfrak{h}$ trivially, there exists a maximal ideal given by    
\begin{align*}
R = R \cap \widetilde{N}^{+} \oplus R \cap \widetilde{N}^{-}.      
\end{align*}
We shall now define the contragredient Lie superalgebra.
  
(2.3) $\mathfrak{g}(A, \tau) = \widetilde{\mathfrak{g}}(A, \tau)/R$.\\      
The images of $e_i, \ f_i$ and $\mathfrak{h}$ will be denoted by the same symbols. We call $e_i$ and $f_i$ as the Chevalley generators. $A$ is called the Cartan matrix and $\mathfrak{h}$ is called the Cartan subalgebra of 
$\mathfrak{g}(A, \tau)$.
     
(2.4) We have a root space decomposition of $\mathfrak{g}(A, \tau)$ given by
\begin{align*}
\mathfrak{g}(A, \tau) = \bigoplus_{\alpha \in \mathfrak{h}^*}\mathfrak{g}(A, \tau)_{\alpha}    
\end{align*}
where $\mathfrak{g}(A, \tau)_{\alpha} = \{ x \in \mathfrak{g}(A, \tau) \ | \ [h,x]= \alpha(h)x \ \forall \ h \in \hh \}$.\\          
Since $\alpha_1, \cdots, \alpha_l$ are linearly independent, we see that 
$\mathfrak{g}(A, \tau)_0 = \mathfrak{h}$.\\   
Let $\Delta = \{ \alpha \in \hh^*\ |\ \alpha \neq 0, \ \mathfrak{g}(A, \tau)_{\alpha} \neq (0)\}$. We shall refer to the elements of $\Delta$ as roots of $\mathfrak{g}(A, \tau)$.

(2.5) $\alpha_1, \cdots, \alpha_l$ are called the simple roots of
$\mathfrak{g}(A, \tau)$. The following are standard facts related to
$\mathfrak{g}(A, \tau)$.       
\begin{enumerate}
\item dim $\mathfrak{g}(A, \tau)_{\alpha} < \infty$.    
\item Any $\alpha \in \Delta$ can be written as a non-negative linear combination of $\alpha_1, \cdots, \alpha_l$ or non-positive linear combination of $\alpha_1, \cdots, \alpha_l$. Thus we have $\Delta= \Delta^{+} \cup \Delta^{-}$ (disjoint union) where $\Delta^{+}$ (respectively $\Delta^{-}$) denotes the set of all positive (respectively negative) roots.  
\item Each root space is homogeneous with $\mathbb{Z}_2$-grading given by\\
$\Delta_{\overline{0}}= \{\alpha \in \Delta \ | \ \mathfrak{g}(A, \tau)_{\alpha} \subseteq \mathfrak{g}(A, \tau)_{\overline{0}} \}$ and\\ $\Delta_{\overline{1}}= \{\alpha \in \Delta \ | \ \mathfrak{g}(A, \tau)_{\alpha} \subseteq \mathfrak{g}(A, \tau)_{\overline{1}} \}$. The roots of $\Delta_{\overline{0}}$ are called even and denoted by $|\alpha|=\overline{0}$ if $\alpha \in \Delta_{\overline{0}}$. Similarly $|\alpha|=\overline{1}$ if $\alpha \in \Delta_{\overline{1}}$.
\item $\mathfrak{g}(A, \tau) = N^{+} \oplus \mathfrak{h} \oplus N^{-}$ where $N^{\pm}= \bigoplus_{\alpha \in \Delta^{\pm}}\mathfrak{g}(A, \tau)_{\alpha}$.                                           
\end{enumerate}

(2.6) There exists an automorphism of $\mathfrak{g}(A, \tau)$, denoted by $\omega$, satisfying\\ 
$\omega(e_i)=-f_i, \ \omega(f_i)=-(-1)^{|\alpha_i|}e_i$ and $\omega(h)=h \ \forall \ h \in \mathfrak{h}$.\\ 
This automorphism is of order 4 if $\tau \neq \phi$ and is of order 2 if $\tau= \phi$. Note that $\omega^2$=Id on $\mathfrak{g}(A, \tau)_{\overline{0}}$ and $\omega^2$=$-$Id on $\mathfrak{g}(A, \tau)_{\overline{1}}$.
     
(2.7) A Cartan matrix $A$ is said to be symmetrizable if there exists a diagonal matrix $D$ such that $AD$ is symmetric (the diagonal matrix consists of non-zero entries).\\
Throughout the paper, we assume that the Cartan matrix is symmetrizable.
      
(2.8) Let us suppose that $\mathfrak{g}(A, \tau)$ is a contragredient Lie superalgebra with symmetrizable matrix $A$. Then there exists a unique (up to a scalar multiple) non-degenerate bilinear form $(\cdot , \cdot)$ on $\mathfrak{g}(A, \tau)$ (see Theorem 5.4.1 together with Remark 5.4.2 and Lemma 8.3.1 of \cite{M}) satisfying the following.                     
\begin{enumerate}
\item $(\cdot , \cdot)$ is non-degenerate when restricted to $\mathfrak{h}$.        
\item Invariant: $([x,y],z)=(x,[y,z]) \ \forall \ x, y, z \in \mathfrak{g}(A, \tau)$.    
\item Super-symmetric: $(x,y)=(-1)^{|x||y|}(y,x)$ for homogeneous elements $x, y$ in $\mathfrak{g}(A, \tau)$.  
\item $(\mathfrak{g}_{\alpha}, \mathfrak{g}_{\beta}) \neq 0 \implies \ \alpha + \beta =0$ and the form restricted to $(\mathfrak{g}_{\alpha}, \mathfrak{g}_{-\alpha})$ is non-degenerate.  
\item Even: $(\mathfrak{g}_{\overline{0}}, \mathfrak{g}_{\overline{1}})=0= (\mathfrak{g}_{\overline{1}}, \mathfrak{g}_{\overline{0}})$.             
\end{enumerate}

(2.9) Since the form is non-degenerate on $\mathfrak{h}$, then for each $\alpha \in \mathfrak{h}^{*}$, there exists a unique element $h_{\alpha} \in \mathfrak{h}$ such that $(h_{\alpha},h)= \alpha(h) \ \forall \ h \in \mathfrak{h}$.\\ 
Consequently the map $\alpha \mapsto h_{\alpha}$ gives rise to a linear isomorphism from $\mathfrak{h}^{*}$ to $\mathfrak{h}$. This clearly induces a non-degenerate bilinear form on $\mathfrak{h}^{*}$ such that $(\alpha,\beta)=(h_{\alpha},h_{\beta}) \ \forall \ \alpha, \beta \in \mathfrak{h}^{*}$.
  
(2.10) Let $\rho \in \mathfrak{h}^{*}$ such that $2(\rho, \alpha_i)=(\alpha_i, \alpha_i) \ \forall \ 1 \leqslant i \leqslant l$. Observe that this $\rho$ is not uniquely determined, but we can fix one such $\rho$. Moreover it is also clear that $h_{\rho} \in \mathfrak{h}$.
                              
(2.11) We shall now define the Casimir operator.\\ 
Let $\mathfrak{g}(A, \tau)$ be a contragredient Lie superalgebra with symmetrizable matrix $A$. Let $(\cdot, \cdot)$ be a bilinear form on 
$\mathfrak{g}(A, \tau)$ satisfying (2.8). Let $\{h_1, \cdots, h_k \}$ be a basis and $\{h^{1}, \cdots, h^{k} \}$ be its dual basis in the sense that $(h_i,h^j)= \delta_{ij} \ \forall \ 1 \leqslant i,j \leqslant k$. Let $\alpha \in \Delta^{+}$ and $k_{\alpha}$=dim$\mathfrak{g}(A, \tau)_{\alpha}$. Now choose a basis $\{e_{\alpha}^{i}\}_{i=1}^{k_{\alpha}}$ of $\mathfrak{g}(A, \tau)_{\alpha}$ and let $\{f_{\alpha}^{i}\}_{i=1}^{k_{\alpha}}$ be a dual basis of $\mathfrak{g}(A, \tau)_{-\alpha}$ in the sense that $(e_{\alpha}^{i},f_{\alpha}^{j})= \delta_{ij}$. Note that $(e_{\alpha}^{i},f_{\alpha}^{i})=1=(f_{\alpha}^{i},e_{\alpha}^{i})$, if $\alpha$ is an even root and $(e_{\alpha}^{i},f_{\alpha}^{i})=1=-(f_{\alpha}^{i},e_{\alpha}^{i})$, if $\alpha$ is an odd root. Define\\
$\Omega_{\alpha}=2\sum_{i=1}^{k_{\alpha}}f_{\alpha}^{i}e_{\alpha}^{i}$ and $\Omega_0=2h_{\rho} + \sum_{i=1}^{k}h_ih^i$.\\
Then define the Casimir operator                                      
\begin{align*}
\Omega= \Omega_0 + \sum_{\alpha \in \Delta^{+}} \Omega_{\alpha}.    
\end{align*}
Note that the definition of $\Omega$ involves an infinite sum and hence cannot be an element of the UEA of $\mathfrak{g}(A, \tau)$. But one can define the completion of $U(\mathfrak{g})$ which we shall denote by $\overline{U}(\mathfrak{g})$ (see \cite{K} for details on $\overline{U}(\mathfrak{g})$ where it is done for Lie algebras, but this can be also extended for Lie superalgebras).

(2.12) $V$ is said to be a weight module for $\mathfrak{g}(A, \tau)$ if
\begin{align*}
\displaystyle{V = \bigoplus_{\lambda \in \mathfrak{h}^* } {V_{\lambda}}}
\end{align*}  
where $V_{\lambda} = \{v \in V |\,\, h.v = \lambda(h)v \,\, \forall \,\, h \in \mathfrak{h} \}$ and $\mathrm{dim}V_{\lambda} < \infty$.
                  
(2.13) Category $\mathcal{O}$\\    
 A weight module $V$ for $\mathfrak{g}(A, \tau)$ is said to be in category $\mathcal{O}$ if for any $v \in V$, $\mathfrak{g}(A, \tau)_{\alpha}.v=0$ for all $\alpha \in \Delta_{+}$ with ht$(\alpha) >>0$.  
\bthm \label{T2.1}
If $V \in \mathcal{O}$, then $(\Omega X - X \Omega).V=0 \ \forall \ X \in \mathfrak{g}$.       
\ethm
The proof of this theorem will be postponed to the next section where we shall prove a much more general result.
 
(2.14) For the rest of this section, we shall take $\mathfrak{g}$ to be a basic classical Lie superalgebra. For simplicity, let us assume that $\mathfrak{g}$ is not of type $A(m,n)$. They are all known to be symmetrizable contragredient Lie superalgebras and hence admit non-degenerate bilinear forms satisfying (2.8). Let\\
$\mathfrak{g} = \mathfrak{h} \bigoplus \big ( \sum_{\alpha \in \Delta} \mathfrak{g}_{\alpha} \big )$ be a root space decomposition of $\mathfrak{g}$. We first give a definite expression to $\rho$ (see (2.9) for the definition of $\rho$). Let $\Delta^{+}= \Delta_{\overline{0}}^{+} \cup \Delta_{\overline{1}}^{+}$ where $\Delta_{\overline{0}}^{+}$ are the even positive roots and $\Delta_{\overline{1}}^{+}$ are the odd positive roots. Define $\rho_{\overline{0}}= \dfrac{1}{2} \sum_{\alpha \in \Delta_{\overline{0}}^{+}} \alpha$ and $\rho_{\overline{1}}= \dfrac{1}{2} \sum_{\alpha \in \Delta_{\overline{1}}^{+}} \alpha$ together with $\widetilde{\rho}= \rho_{\overline{0}} - \rho_{\overline{1}}$.  
\bppsn\cite[Proposition 1.33]{CW} \label{P2.2} 
$2(\widetilde{\rho}, \alpha_i)=(\alpha_i, \alpha_i) \ \forall \ i \in I$. 
\eppsn                    
This proposition allows us to take $\rho= \widetilde{\rho}$. We now recall a result from \cite{M}.
\bppsn\cite[Lemma 8.5.1]{M}\label{P2.3} 
Let $\mathfrak{g}$ be a Lie superalgebra as in (2.14) and $(\cdot, \cdot)$ be a non-degenerate bilinear form satisfying (2.8). Let $\{x_1, \cdots, x_k \}$ be a basis of $\mathfrak{g}$ and $\{y_1, \cdots, y_k \}$ be the corrresponding dual basis (where we take dim $\mathfrak{g}=k$) with $(x_i,y_j)= \delta_{ij}$ and $x_i, \ y_i$ are homogeneous elements of degree $\beta_i$. If we now define
$\Omega_c = \sum_{i=1}^{k} (-1)^{\beta_i}x_i y_i \in U(\mathfrak{g})$, then $\Omega_c$ commutes with $U(\mathfrak{g})$.          
\eppsn

(2.15) We shall prove that $\Omega_c= \Omega$ where $\Omega$ is defined in (2.11). Let $\{h_1, \cdots, h_m \}$ be a basis of $\mathfrak{h}$ and $\{h^1, \cdots, h^m \}$ be its dual basis. For $\alpha \in \Delta^{+}$, let $\{e_{\alpha}^{i} \}$ be a basis of $\mathfrak{g}_{\alpha}$ and $\{f_{\alpha}^{i} \}$ be the dual basis of $\mathfrak{g}_{-\alpha}$ in the sense that $(e_{\alpha}^{i}, f_{\alpha}^{j})= \delta_{ij}$. Observe that $(e_{\alpha}^{i}, f_{\alpha}^{i})=1=(f_{\alpha}^{i}, e_{\alpha}^{i})$ for even root $\alpha$ and $(e_{\alpha}^{i}, f_{\alpha}^{i})=1=-(f_{\alpha}^{i}, e_{\alpha}^{i})$ for odd root $\alpha$. Finally we shall also consider the basis $\{h_1, \cdots, h_m \} \cup \{e_{\alpha}^{i} \}_{\alpha \in \Delta^{+}} \cup \{f_{\alpha}^{i} \}_{\alpha \in \Delta^{+}}$ of $\mathfrak{g}$ and take its dual basis as $\{h^1, \cdots, h^m \} \cup \{f_{\alpha}^{i} \}_{\alpha \in \Delta^{+}} \cup \{(-1)^{|\alpha|}e_{\alpha}^{i} \}_{\alpha \in \Delta^{+}}$.        
\bppsn
$\Omega_c= \Omega$ where the notations are as described above.        
\eppsn
\begin{proof}
First recall that $[x_{\alpha},y_{\alpha}]=(x_{\alpha},y_{\alpha})h_{\alpha}$ from \cite[Lemma 8.3.2]{M} where $x_{\alpha} \in \mathfrak{g}_{\alpha}, \ y_{\alpha} \in \mathfrak{g}_{-\alpha}$ and $h_{\alpha}$ is as defined in (2.9). Then we have 
\begin{align*}
\Omega_c= \sum_{i}h_ih^i + \sum_{\alpha \in \Delta_{\overline{0}}^{+}}e_{\alpha}^{i}f_{\alpha}^{i} + \sum_{\alpha \in \Delta_{\overline{0}}^{+}}f_{\alpha}^{i}e_{\alpha}^{i} + \sum_{\alpha \in \Delta_{\overline{1}}^{+}}(-1)e_{\alpha}^{i}f_{\alpha}^{i} \\ + \sum_{\alpha \in \Delta_{\overline{1}}^{+}}(-1)f_{\alpha}^{i}(-1)e_{\alpha}^{i}.  
\end{align*}
Note that $e_{\alpha}^{i}f_{\alpha}^{i} = f_{\alpha}^{i}e_{\alpha}^{i} + h_{\alpha}$ if $\alpha \in \Delta_{\overline{0}}^{+}$ and $e_{\alpha}^{i}f_{\alpha}^{i} =-f_{\alpha}^{i}e_{\alpha}^{i} + h_{\alpha}$ if $\alpha \in \Delta_{\overline{1}}^{+}$.\\
This directly implies that
\begin{align*}
\Omega_c= \sum_{i}h_ih^i + \sum_{\alpha \in \Delta_{\overline{0}}^{+}}2f_{\alpha}^{i}e_{\alpha}^{i} + \sum_{\alpha \in \Delta_{\overline{0}}^{+}}h_{\alpha} + \sum_{\alpha \in \Delta_{\overline{1}}^{+}}2f_{\alpha}^{i}e_{\alpha}^{i} - \sum_{\alpha \in \Delta_{\overline{1}}^{+}}h_{\alpha}.               
\end{align*}
To complete the proof, we need to show that $2h_{\rho}= \sum_{\alpha \in \Delta_{\overline{0}}^{+}}h_{\alpha} - \sum_{\alpha \in \Delta_{\overline{1}}^{+}}h_{\alpha}$. But this follows from Proposition \ref{P2.2} and the definition of $h_{\alpha}$ in (2.9).        
\end{proof}
The following lemma will be useful for calculation.
\blmma \label{L2.5}
Let $\mathfrak{g}$ be any Lie superalgebra and suppose that $x,y,z$ are homogeneous elements of $\mathfrak{g}$. Then
\begin{enumerate}
\item $[z,xy]=[z,x]y + (-1)^{|z||x|}x[z,y]$.
\item $[xy,z]=x[y,z] +  (-1)^{|y||z|}[x,z]y$.   
\end{enumerate}      
\elmma
(2.16) We shall now work out $\Omega$ for rank 1 Lie superalgebras which will be utilized in the next section to prove that the Casimir operator is central.\\
(1) $\mathfrak{g}=\mathfrak{sl}_{2}-$We know that $\mathfrak{sl}_{2}$ is spanned by $\{e,h,f\}$ satisfying the relations\\
$[e,f]=h, \ [h,e]=2e$ and $[h,f]=-2f$.\\
The symmetric bilinear form on $\mathfrak{sl}_{2}$ is given by\\
$(e,f)=0, \ (h,h)=2$ and $(e,h)=0=(f,h)$.\\
In this case, the Casimir element is given by
$\Omega_{\mathfrak{sl}_{2}} = h +h^2/2 + 2fe$.\\
It is easy to check that $\Omega_{\mathfrak{sl}_{2}}$ is central.\\
(2) $\mathfrak{g}=\mathfrak{gl}(1,1)-$This Lie superalgebra is not covered earlier as it is solvable. But $\mathfrak{g}$ is endowed with a non-degenerate bilinear form satisfying (2.8). The Casimir operator can be defined as before which can be written explicitly.\\
Let 
\[e=
\begin{bmatrix}
0 & 1 \\
0 & 0
\end{bmatrix},               
 \ f=
\begin{bmatrix}
0 & 0 \\
1 & 0
\end{bmatrix}
\ \text{and} \ \ h_1=
\begin{bmatrix}
1 & 0 \\
0 & 0
\end{bmatrix},
 \ h_2=
\begin{bmatrix}
0 & 0 \\
0 & 1
\end{bmatrix}.               
\]        
Set $h=h_1 + h_2$. Then $\mathfrak{g}$ is spanned by the elements $\{e,f,h_1,h_2\}$ that satisfy the following relations.\\
$[e,f]=h=[f,e], \ [e,h]=0=[f,h]$ along with\\
$-[e,h_1]=e=[e,h_2], \ [f,h_1]=f=-[f,h_2]$.\\
Observe that there is only one root which is isotropic and the non-degenerate bilinear form on $\mathfrak{g}$ is given by $(e,f)=1=-(f,e), \ (h_1,h_1)=1=-(h_2,h_2),\\ (h_1,h_2)=0$ with all the other values being zero. It is also clear that\\
$\Omega_{\mathfrak{gl}(1,1)} = -h + h_{1}^{2} - h_{2}^{2} + 2fe$ which coincides with the earlier definition given for the other Lie superalgebras. Now from Lemma \ref{L2.5}, it follows that $[e,fe]=[e,f]e$ and $[f,fe]=-f[f,e]$, from which we can easily conclude that $\Omega_{\mathfrak{gl}(1,1)}$ is central in $U(\mathfrak{g})$.\\
(3) $\mathfrak{g}=\mathfrak{osp}(1,2)$ (Lie superalgebra of type $B(0,1)$)$-$ This is a contragredient Lie superalgebra admits of a non-degenerate bilinear form satisfying (2.8). This Lie superalgebra have already been covered earlier and we have also previously defined the Casimir operator $\Omega= \Omega_{\mathfrak{osp}(1,2)}$ in (2.11).\\
The Lie superalgebra $\mathfrak{g}$ is of rank one having Chevalley generators $e$ and $f$ satisfying $[e,f]=h$. Let $\alpha$ be the root for the root vector $e$ where $\alpha$ is odd and non-isotropic. Then $2\alpha$ is an even root. Putting $4e^{\prime}=[e,e]$ and $-4f^{\prime}=[f,f]$, it is now evident that $\mathfrak{g}$ is spanned by $\{e,f,e^{\prime},f^{\prime},h \}$. In this instance, the super Lie bracket is given by the following table which is reproduced from \cite{W}. Here the brackets are computed by simply taking the first entry from the first column and then taking the other entry (which is now fixed) from the first row.         
\[
\begin{array}{l|*{5}{l}}
& \ \ h   & \ \ e^{\prime}   & \ \ f^{\prime} &  \ \ e  & \ \ f \\
\hline
h            & \ \ 0    & \ \ 4e^{\prime}   & -4f^{\prime} & \ \ 2e  & -2f \\
e^{\prime}   & -4e^{\prime}   & \ \ 0 & \ h/2 & \ \ 0  & -e \\
f^{\prime}   & \ \ 4f^{\prime} & -h/2 & \ \ 0 & -f  & \ \ 0 \\
e            & -2e & \ \ 0 & \ \ f & \ \ 4e^{\prime}  & \ \ h \\ 
f            & -2f & \ \ e & \ \ 0 & \ \ h  & -4f^{\prime} \\
\end{array} 
\]
In this case, $e,f$ are odd vectors and $e^{\prime},f^{\prime},h$ are even vectors.\\
The non-degenerate bilinear form is given by
$(e,f)=1=-(f,e), \ (h,h)=2$ and $(e^{\prime},f^{\prime})=1/4=(f^{\prime},e^{\prime})$. The remaining values are zero.\\
It is easy to see that $h_{\rho}=2h-h=h$ and the Casimir operator is given by\\
$\Omega= \Omega_{\mathfrak{osp}(1,2)} = h + h^2/2 + 8f^{\prime}e^{\prime} + 2fe$. Now from Proposition \ref{P2.2} and (2.15), we see that $\Omega= \Omega_c$ which also commutes with $U(\mathfrak{g})$. One can verify that $\Omega$ is central by using the following identities that follow from Lemma \ref{L2.5}.\\
$[e,fe] = [e,f]e - f[e,e], \  [f,fe] = [f,f]e - f[f,e]$ and\\
$[e,f^{\prime}e^{\prime}] = [e,f^{\prime}]e^{\prime} + f^{\prime}[e,e^{\prime}] = [e,f^{\prime}]e^{\prime}, \ [f,f^{\prime}e] = [f,f^{\prime}]e + f^{\prime}[f,e] = f^{\prime}[f,e]$.\\
Here $\{h,e,f,e^{\prime},f^{\prime} \}$ is a basis of $\mathfrak{g}$ and $\{h/2,f,-e,4f^{\prime},4e^{\prime} \}$ is its dual basis.           
        
\section{Casimir Operators and Loop Algebras}\label{S3}

In this section, we define the notion of loop algebra of a contragredient Lie superalgebra and then introduce the notion of central operators which we call generalized Casimir operators.    

(3.1) Suppose that $\mathfrak{g}= \mathfrak{g}(A, \tau)$ is a symmetrizable contragredient Lie superalgebra equipped with a Cartan subalgebra $\mathfrak{h}$ and let $(\cdot, \cdot)$ be a non-degenerate bilinear form satisfying (2.8). Let $\{\alpha_1, \cdots, \alpha_l \}$ be a set of simple roots and $\{\alpha_{1}^{\vee}, \cdots, \alpha_{l}^{\vee} \}$ be its co-roots. Let $\Delta$ and $\Delta_{+}$ denote the respective collection of all roots and positive roots of $\mathfrak{g}$ with respect to $\mathfrak{h}$. We also have the corresponding root space decomposition given by $\mathfrak{g} = \mathfrak{h} \bigoplus \sum_{\alpha \in \Delta} \mathfrak{g}_{\alpha}$.  

(3.2) Let $B$ be an associative, commutative and finitely generated algebra with unity. For any vector space $V$, put $V(B)=V \otimes B$. Then $\mathfrak{g}(B)$ carries a natural Lie superalgebra structure with even part $\mathfrak{g}_{\overline{0}}(B)$ and odd part $\mathfrak{g}_{\overline{1}}(B)$. Let $\mathfrak{g}= N^{-} \oplus \mathfrak{h} \oplus N^{+}$ be the standard triangular decomposition. Then $\mathfrak{g}(B)= N^{-}(B) \oplus \mathfrak{h}(B) \oplus N^{+}(B)$ is a triangular decomposition of $\mathfrak{g}(B)$. For $\alpha = \sum_{i=1}^{l}  m_i \alpha_i \in \Delta_{+}$, let ht($\alpha$)= $\sum_{i=1}^{l} m_i$. Also note that $\mathfrak{g} \cong \mathfrak{g} \otimes \mathbb{C}$.\\
We shall now define a category of modules for the loop algebra $\mathfrak{g}(B)$, similar to the one defined earlier in (2.13) and denote this category again by $\mathcal{O}$.    
 
(3.3) A module $V$ over $\mathfrak{g}(B)$ is said to be in the category $\mathcal{O}$ if the following conditions are satisfied.
\begin{enumerate}
\item $V$ is a weight module with respect to $\mathfrak{h}$ (see (2.12)). \item For any $v \in V$ and $a \in B$, we have $(x_{\alpha} \otimes a).v=0$ for ht($\alpha$)$>>0$ where $x_{\alpha} \in \mathfrak{g}_{\alpha}, \alpha \in \Delta_{+}$.          
\end{enumerate}

(3.4) We shall now extract a class of irreducible modules over $\mathfrak{g}(B)$ in the category $\mathcal{O}$. Consider the $1$-dimensional $\big (\mathfrak{h}(B) \oplus N^{+}(B) \big )$-module $\mathbb{C}v$ where $N^{+}(B)$ acts trivially on this module and $\mathfrak{h}(B)$ acts by $\psi \in \mathfrak{h}(B)^*$ on $\mathbb{C}v$. Let
\begin{align*}
M(\psi)=U(\mathfrak{g}(B)) \bigotimes_{\mathfrak{h}(B) \oplus N^{+}(B)} \mathbb{C}v
\end{align*}    
be the corresponding Verma module. By standard arguments, it follows that $M(\psi)$ has a unique maximal submodule and $M(\psi)$ has a unique irreducible quotient which we shall denote by $V(\psi)$. Note that when $B$ is infinite dimensional, $M(\psi)$ does not have finite dimensional weight spaces. $V(\psi)$ has finite dimensional weight spaces depending on $\psi$.

(3.5) Set $\mathfrak{g}^{\prime}=[\mathfrak{g}, \mathfrak{g}]$ and $\mathfrak{h}^{\prime}= \mathfrak{g}^{\prime} \cap \mathfrak{h}$. Let
$\mathfrak{h}^{\prime \prime}$ be any vector subspace of $\mathfrak{h}$ such that $\mathfrak{h} = \mathfrak{h}^{\prime} \oplus \mathfrak{h}^{\prime \prime}$. Put $\widetilde{\mathfrak{g}} = \mathfrak{g}^{\prime}(B) \oplus \mathfrak{h}^{\prime \prime}$. It is not too difficult to see that $V(\psi)$ is an irreducible module for $\mathfrak{g}^{\prime}(B)$ (see \cite[Lemma 1.5]{R2}).
   
\bppsn\cite[Proposition 1.6]{R2}  
$V(\psi)$ has finite dimensional weight spaces if and only if there exists a
co-finite $I$ of $B$ (i.e. dim($B/I$) $<\infty$) such that $\psi({\mathfrak{h}}^{\prime} \otimes I)=0$. In this case, we have $\mathfrak{g}^{\prime}(I).V(\psi)=(0)$. 
\eppsn

(3.6) In the above instance, $V(\psi) \in \mathcal{O}$.

(3.7) Refer to (1.8) in \cite{R2} for some special types of co-finite ideals. 

(3.8) \textbf{Central Operator.} A linear map $T$ on an object of $\mathcal{O}$ is called a central operator if $T$ commutes with the action of $\mathfrak{g}$.\\
We shall now define a class of central operators $\Omega(a,b)$ where $a,b \in B$ and call them generalized Casimir operators. These operators are similar to the ones defined in \cite{R2} for the Lie algebra case.\\
For $\alpha \in \Delta_{+}$, let $\{e_{\alpha}^{i} \}$ be a basis of $\mathfrak{g}_{\alpha}$ and $\{f_{\alpha}^{i} \}$ be the dual basis of $\mathfrak{g}_{-\alpha}$ in the sense that $(e_{\alpha}^{i}, f_{\alpha}^{j})= \delta_{ij}$. Note that $(e_{\alpha}^{i}, f_{\alpha}^{i})=1=(f_{\alpha}^{i}, e_{\alpha}^{i})$ for even root $\alpha$ and $(e_{\alpha}^{i}, f_{\alpha}^{i})=1=-(f_{\alpha}^{i}, e_{\alpha}^{i})$ for odd root $\alpha$. Let $\{h_i \}$ be a basis of $\mathfrak{h}$ and $\{h^i \}$ be the dual basis of $\mathfrak{h}$ such that $(h_i,h^j)= \delta_{ij}$. Let us further set $x(a)= x \otimes a$ for $x \in \mathfrak{g}$ and $a \in B$.\\
Define $\Omega_{\alpha}(a,b) = \sum_{i} f_{\alpha}^{i}(a)e_{\alpha}^{i}(b)$ and $\Omega(a,b) = 2h_{\rho}(ab) + \sum_{i}h_i(a)h^i(b) + \sum_{\alpha \in \Delta^{+}}\Omega_{\alpha}(a,b) + \sum_{\alpha \in \Delta^{+}}\Omega_{\alpha}(b,a)$ where $h_{\rho}$ is defined in (2.10). Observe that $\Omega(a,b)=\Omega(b,a)$ and $\Omega(1,1)= \Omega$ (as defined in the previous section).            

\bthm \label{T3.2}
$\Omega(a,b)$ is central for any $a,b \in B$.
\ethm
The proof will be provided towards the end of this section as we need some preparation to prove this theorem. 

\brmk
The proof of Theorem \ref{T2.1} follows from the above theorem by just taking 
$B= \mathbb{C}$ and $a=b=1$.
\ermk 

(3.9) Let us first verify Theorem \ref{T3.2} for the rank $1$ contragredient Lie supralgebras. There are three of them, namely $\mathfrak{sl}_2, \ \mathfrak{gl}(1,1)$ and $\mathfrak{osp}(1,2)$. We shall follow the notations as in (2.16).\\
(1) $\mathfrak{sl}_2-$ $\Omega_{\mathfrak{sl}_{2}}(a,b) = h(ab) + \dfrac{1}{2}h(a)h(b) + f(a)e(b) + f(b)e(a)$.\\
It can be readily checked that $\Omega_{\mathfrak{sl}_{2}}(a,b)$ is central.\\
(2) $\mathfrak{gl}(1,1)-$ $\Omega_{\mathfrak{gl}(1,1)}(a,b) = -h(ab) + h_1(a)h_1(b) - h_2(a)h_2(b) + f(a)e(b) + f(b)e(a)$.
It is easy to see that $\Omega_{\mathfrak{gl}(1,1)}$ is central.\\
(3) $\mathfrak{osp}(1,2)-$ $\Omega_{\mathfrak{osp}(1,2)} = h(ab) + \dfrac{1}{2}h(a)h(b) + 4f^{\prime}(a)e^{\prime}(b) + 4f^{\prime}(b)e^{\prime}(a) + f(a)e(b) + f(b)e(a)$. It can be verified that this operator is central.    

Next let us prove the following proposition whose proof is quite similar to \cite[Lemma 18.4.1]{M}, but with a sign difference.    
            
(3.10) For $\alpha, \beta \in \Delta_{+}$, let $\{x_i \}$ and $\{u_i \}$ be
bases of $\mathfrak{g}_{\alpha}$ and $\mathfrak{g}_{\beta}$ respectively. Further suppose that $\{y_i \}$ and $\{v_i \}$ are bases of $\mathfrak{g}_{-\alpha}$ and $\mathfrak{g}_{-\beta}$ respectively satisfying $(x_i,y_i)= \delta_{ij}$ and $(u_i,v_j)= \delta_{ij}$.       

\bppsn \label{P3.4}
Let $z \in \mathfrak{g}_{\beta - \alpha}$ with $|z|=|\beta - \alpha|$. Then
$\sum_{i} [v_i,z] \otimes u_i = (-1)^{|z|} \sum_{i} y_i \otimes [z,x_i]$.
\eppsn
\begin{proof}
First note that both sides belong to $\mathfrak{g}_{- \alpha} \otimes \mathfrak{g}_{\beta}$. Define a bilinear form on $\mathfrak{g} \otimes \mathfrak{g}$ by setting $(u \otimes v, x \otimes y) = (u,x)(v,y)$. It is clear that $(\cdot, \cdot)$ induces a non-degenerate pairing between $\mathfrak{g}_{\alpha} \otimes \mathfrak{g}_{\beta}$ and $\mathfrak{g}_{{\alpha}^{\prime}} \otimes \mathfrak{g}_{{\beta}^{\prime}}$ where $\alpha + \alpha^{\prime} = 0 = \beta + \beta^{\prime}$.\\
In order to prove this proposition, it suffices to check that
\begin{align*}
\big (\sum_{i} [v_i,z] \otimes u_i - (-1)^{|z|} \sum_{i} y_i \otimes [z,x_i], a \otimes b \big ) = 0 \ \forall \ a \in \mathfrak{g}_{\alpha}, \ b \in \mathfrak{g}_{-\beta}.   
\end{align*}
Now using the orthogonality property, we have $b = \sum_{i}(u_i,b)v_i$ which implies that $[b,z] = \sum_{i}(u_i,b)[v_i,z]$. This in turn immediately gives us\\
$([b,z],a) = \sum_{i}([v_i,z],a)(u_i,b) = \sum_{i}([v_i,z] \otimes u_i, a \otimes b) \cdots \cdots$ ($i$)\\
Again by the orthogonality property, we can write
$[z,b] = \sum_{i}(x_i,[z,b])y_i$.\\
Henceforth by applying the invariance property of the form on $\mathfrak{g}$, we obtain\\
$([z,b],a) = \sum_{i}(x_i,[z,b])(y_i,a) = \sum_{i}([x_i,z],b)(y_i,a)= \sum_{i}(y_i,a)([x_i,z],b)$\\ 
$ = \sum_{i} ( y_i \otimes [x_i,z], a \otimes b) = \sum_{i} (-1)(-1)^{|\alpha||z|} ( y_i \otimes [z,x_i], a \otimes b) \cdots \cdots$ ($ii$)\\
But since $([b,z],a) = (-1)(-1)^{|z||\beta|}([z,b],a)$ and $|\alpha - \beta| |\alpha - \beta| = |\alpha - \beta| = |z|$, we finally get (from ($ii$)) $([b,z],a) = \sum_{i} (-1)^{|z|} ( y_i \otimes [z,x_i], a \otimes b) \cdots \cdots$ ($iii$)\\
The proposition is now a direct consequence of ($i$) and ($iii$).     
\end{proof}

\bcrlre \label{C3.5}
For $a, b \in B, \ \sum_{i}[v_i(a),z]u_i(b) = (-1)^{|z|} \sum_{i} y_i(a)
[z,x_i(b)]$.
\ecrlre
\begin{proof}
Consider the $\mathfrak{g}$-module homomorphism given by
\begin{align*}
{\phi} : \mathfrak{g} \otimes \mathfrak{g} \longrightarrow U(\mathfrak{g} \otimes A)
\end{align*}
$\,\,\,\,\,\,\,\,\,\,\,\,\,\,\,\,\,\,\,\,\,\,\,\,\,\,\,\,\,\,\,\,\,\,\,\,\,\,\,\,\,\,\,\,\,\,\,\,\,\,\,\,\,\,\,\,\,\,\,\,\,\,\,\,\,\,\,\,\,\,\,\,\,\,\,\,\ x \otimes y \longmapsto x(a)y(b)$\\
Applying this map in Proposition \ref{P3.4}, we thereby obtain the desired result.     
\end{proof}
\textbf{Proof of Theorem \ref{T3.2}.} We shall closely follow the proof presented in \cite[Theorem 18.4.2]{M} although our Casimir operator differs from \cite{M}.\\ For fixed $1 \leqslant j \leqslant l$, set $\Delta_j = \mathbb{Z}_{+}e_j \cap \Delta_{+}$ and $\theta_j = \Delta_{+} \setminus \Delta_j$. We want to prove that $\Omega e_j = e_j \Omega$. In order to accomplish this, let us first show that
\begin{align*}
\sum_{\alpha \in \theta_j} \big (\Omega_{\alpha}(a,b)e_j - e_j \Omega_{\alpha}(a,b) \big ) = 0.
\end{align*} 
As $\Omega_{\alpha}(a,b)$ is an even vector, the LHS of the above equation is equal to $\sum_{\alpha \in \Delta_j} [\Omega_{\alpha}(a,b), e_j]$, whence by Lemma \ref{L2.5}(2), this term is again equal to 
   
(3.11) \,\,\,\,\,\,\,\ $\sum_{\alpha \in \theta_j}\big(\sum_{i} (f_{\alpha}^{i}(a) [e_{\alpha}^{i}(b),e_j] + (-1)^{|\alpha||e_j|}[f_{\alpha}^{i}(a),e_j]e_{\alpha}^{i}(b))\big)$.\\
Note that on any vector $v$ in $V \in \mathcal{O}$, the above sum is finite.\\
Take $z=e_j$ and $\beta = \alpha + e_j$. We claim that\\ $(-1)^{|\alpha_j||\alpha + \alpha_j|} \sum_{i} [f_{\beta}^{i}(a),e_j]e_{\beta}^{i}(b) = (-1) \sum_{i} f_{\alpha}^{i}(a)[e_{\alpha}^{i}(b),e_j]$.\\
To prove this claim, we apply Corollary \ref{C3.5} by taking $\alpha, \beta$ with $\alpha - \beta = e_j$.\\
We have $\sum_{i} [f_{\beta}^{i}(a),e_j]e_{\beta}^{i}(b) = (-1)^{|\alpha_j|}\sum_{i} f_{\alpha}^{i}(a)[e_j,e_{\alpha}^{i}(b)]$\\
$= (-1)(-1)^{|\alpha_j| + |\alpha_j||\alpha|} \sum_{i} f_{\alpha}^{i}(a)[e_{\alpha}^{i}(b),e_j]$.\\
The claim now follows by observing that $|\alpha_j||\alpha + \alpha_j| + |\alpha_j| + |\alpha_j||\alpha| = \overline{0}$.\\
As a result of this claim, it is clear that the cross terms also get cancelled and thus the equation in (3.11) finally reduces to zero. So what remins to be shown is the following. If $\alpha = \alpha_j$, then

(3.12) \,\,\,\ $2h_{\rho}(ab) + \sum_{i}h_i(a)h^i(b) + \sum_{\alpha \in \Delta_j} \big (\sum_{i}(f_{\alpha}^{i}(a) e_{\alpha}^{i}(b) + f_{\alpha}^{i}(b)e_{\alpha}^{i}(a)) \big )$\\
commutes with $e_j$ which is essentially a rank $1$ calculation.\\
\textbf{Case 1.} $\alpha$ is a simple even root.\\
In this case $\Delta_j = \{\alpha_j \}$. Let $[e_{\alpha},f_{\alpha}]=h_{\alpha}$ and consider\\
$h_{\alpha}^{\perp} = \{h \in \mathfrak{h} \ | \ (h,h_{\alpha})=0 \}$ which is equal to Ker($\alpha$).\\
Now one can choose $\{h_i \}$ and $\{h^i \}$ in such a way that $\{h_i \}_{i \geqslant 2}$ and $\{h^i \}_{i \geqslant 2}$ is a dual basis for $h_{\alpha}^{\perp}$ with $h_1$ and $h^1$ being multiples of $h_{\alpha}$. It is easy to check that (3.12) commutes with $e_j$.\\
\textbf{Case 2.} $\alpha$ is a non-isotropic odd root.\\
This phenomenon takes place in case of $\mathfrak{osp}(1,2)$ where $[e_{\alpha},f_{\alpha}]=h_{\alpha}$ and $(h_{\alpha},h_{\alpha})=(\alpha, \alpha) \neq 0$. Let $h_{\alpha}^{\perp} = \{h \in \mathfrak{h} \ | \ (h,h_{\alpha})=0 \}$.\\
We can readily verify that (3.12) commutes with $e_j$ as in \textbf{Case 1} (see (3.9)(3) where $h_{\alpha}$ is denoted by $h$).\\
\textbf{Case 3.} $\alpha$ is an odd isotropic root.\\
Here we deal with $\mathfrak{sl}(1,1)$ where $[e_{\alpha},f_{\alpha}]=h_{\alpha}$ and $(h_{\alpha},h_{\alpha})=(\alpha, \alpha)=0$.\\
First observe that $\mathfrak{sl}(1,1) \neq \mathfrak{g}$ as $\mathfrak{sl}(1,1)$ does not carry any non-degenerate bilinear form satisfying (2.8) which thereby implies that dim $\mathfrak{h} \geqslant 2$. Let us now enlarge $\mathfrak{sl}(1,1)$ to $\mathfrak{gl}(1,1)$. Pick $h_0 \in \mathfrak{h}$ such that $(h_0,h_{\alpha})=1$. Further choose $d \in \mathbb{C}$ such that $h_1 = h_0 + dh_{\alpha}$ and ($h_1,h_1$)=$1$. Take $h_2 = h_{\alpha}-h_1$ and check that $(h_2,h_2)=1$. Putting $h=h_{\alpha}$, we have the following.\\
$[e,h]=0=[f,h]$ as $\alpha$ is an isotropic root. Using the invariance, it can be easily verified that $[e,h_1]=-e, \ [f,h_1]=f, \ [e,h_2]=e, \ [f,h_2]=-f$.\\
Let $\widetilde{\mathfrak{h}} = \{h_1, h_2 \}$ and consider $\widetilde{\mathfrak{h}}^{\perp} = \{h \in \mathfrak{h} \ | \ (h,h_1)=0=(h,h_2) \}$.\\
We claim that $\widetilde{\mathfrak{h}}^{\perp} \subset$ Ker$(\alpha)$. To this end, let us consider\\
$0 = (\widetilde{\mathfrak{h}}^{\perp},[e_{\alpha},f_{\alpha}]) = ([\widetilde{\mathfrak{h}}^{\perp},e_{\alpha}],f_{\alpha}) = \alpha(\widetilde{\mathfrak{h}}^{\perp})(e_{\alpha},f_{\alpha})$.\\
Now since $(e_{\alpha},f_{\alpha})=1$, we obtain $\alpha(\widetilde{\mathfrak{h}}^{\perp})=0$. Hence the claim.\\
The equality clearly does not hold as $h_{\alpha} \in$ Ker$(\alpha)$, but $h_{\alpha} \notin \widetilde{\mathfrak{h}}^{\perp}$.\\
Pick a dual basis of $\mathfrak{h}$ such that $\{h_1, h_2, h_i \ | \ i \geqslant 3 \}$ and $\{h_1, -h_2, \  h^i \ | \ i \geqslant 3 \}$.\\
Observe that $\big [\sum_{i \geqslant 3} h_i(a)h^i(b),e_j \big ]=0$. The rest of the operators in (3.11) is precisely $\Omega_{\mathfrak{gl}(1,1)}$ for $\mathfrak{gl}(1,1)$ and commutes with $e_j$ as noted in (3.9)(2).\\
Thus we have proved that $[\Omega(a,b),e_j]=0$ and a similar argument also yields that $[\Omega(a,b),f_j]=0$. Finally as $\Omega(a,b)$ has weight $0$, it can be trivially deduced that $\Omega(a,b)$ commutes with $\mathfrak{h}$. This proves the theorem.                                              
\section{Tensor Product Modules For General Lie Superalgebras}\label{S4}

In this section, we work with general linear Lie superalgebras. We consider tensor product modules and define central operators that commute with general linear Lie superalgebra. These operators take one highest weight vector to a new highest weight vector.

(4.1) Let us fix $\mathfrak{g} = \mathfrak{gl}(M,N)$ where $M, N \in \mathbb{N}$ and $M+N \geqslant 3$.\\
Consider a $\mathbb{Z}_2$-graded vector space $V = V_{\overline{0}} + V_{\overline{1}}$ where dim$V_{\overline{0}}=M$ and dim$V_{\overline{1}}=N$.
This $\mathbb{Z}_2$-gradation on $V$ naturally extends to the vector space End$V = ($End$V)_{\overline{0}} \oplus$ (End$V)_{\overline{1}}$ by setting\\
(End$V)_{\overline{j}} = \{f \in \text{End}(V) \ | \ f(V_{\overline{k}}) \subseteq V_{\overline{k} + \overline{j}} \ \forall \ \overline{k} \in \mathbb{Z}_2 \}$.\\
Then End$V$ becomes a Lie superalgebra under the super bracket operation\\
$[f,g] = f \circ g - (-1)^{\overline{i} \overline{j}}g \circ f \ \forall \ f \in (\text{End}V)_{\overline{i}} \ , \ g \in (\text{End}V)_{\overline{j}}$.\\
Then we clearly have End$V = \mathfrak{gl}(M,N) = \mathfrak{g}$. By fixing a basis of $V_{\overline{0}}$ and $V_{\overline{1}}$, we can write $f \in$ End$V$ as
\[f=
\begin{bmatrix}
A & B \\
C & D
\end{bmatrix}
\]         
where $A$ is a $M \times M$ matrix, $B$ is a $N \times M$ matrix, $C$ is a $M \times N$ matrix and $D$ is a $N \times N$ matrix.\\
It is easy to see that 
\[
\begin{bmatrix}
A & 0 \\
0 & D
\end{bmatrix} \in (\text{End}V)_{\overline{0}} \ \text{and} 
\begin{bmatrix}
0 & B \\
C & 0
\end{bmatrix} \in (\text{End}V)_{\overline{1}}.  
\]
Let $\mathfrak{h}$ be a Cartan subalgebra of $\mathfrak{g}$ which is contained in (End$V)_{\overline{0}}$ and we take them as diagonal matrices. Then we have a root space decomposition   
\begin{align*}
\mathfrak{g} =  \mathfrak{h} \oplus \sum_{\alpha \in \Delta} \mathfrak{g}_{\alpha}.
\end{align*}
Let $\Delta^{+}$ denote the set of all positive roots with respect to this Cartan subalgebra $\mathfrak{h}$ and set $\mathfrak{g}^{+} = \bigoplus_{\alpha \in \Delta^{+}} \mathfrak{g}_{\alpha}$.    

(4.2) Let $\{E_{ij} \ | \ 1 \leqslant i,j \leqslant M+N \}$ be the standard matrix elements. It is well-known that $E_{ij}$ is a root vector of $\alpha_{ij}$. If $\alpha_{ij}$ is even, we write $|\alpha_{ij}| = \overline{0}$ and if $\alpha_{ij}$ is odd, we write $|\alpha_{ij}| = \overline{1}$. Put I $= \{1,2, \cdots, M+N \}$ along with I$_0= \{1,2, \cdots, M \}$ and I$_1= \{M+1, M+2, \cdots, M+N \}$.\\
Define parity for I as $p(i) = \overline{0}$ if $i \in $ I$_0$ and  $p(i) = \overline{1}$ if $i \in $ I$_1$ so that $|\alpha_{ij}| = p(i) + p(j)$.     

(4.3) If $V_1, \cdots, V_k$ be modules over $\mathfrak{g}$, then we can put a $\mathfrak{g}$-module on $W = \bigotimes_{i=1}^{k}V_i$ in the following way.\\
Let $x \in \mathfrak{g}$ and $w_j \in V_j$ be homogeneous elements of degree $|x|$ and $|w_j|$. Set
\begin{align*}
x.(w_1 \otimes \cdots \otimes w_k) = x.w_1 \otimes w_2 \otimes \cdots \otimes w_k + (-1)^{|x||w|}w_1 \otimes x.w_2 \otimes \cdots \otimes w_k \\ + \cdots \cdots + (-1)^{|x|(|w_1| + \cdots + |w_{k-1}|)}w_1 \otimes \cdots \otimes x.w_k \ (\text{see (3.2.1) in \cite{CW}}).      
\end{align*}

(4.4) We shall now recall certain standard algebras corresponding to $\mathfrak{g}$. Let $T(\mathfrak{g}) = \bigoplus_{k \in \mathbb{Z}_{+}}T^k(\mathfrak{g})$ be the tensor algebra and $U(\mathfrak{g})$ be the universal enveloping algebra. Then $T(\mathfrak{g})$ and $U(\mathfrak{g})$ are $\mathbb{Z}_2$-graded and there exists an onto algebra homomorphism   
\begin{align*}
\pi : T(\mathfrak{g}) \longrightarrow U(\mathfrak{g}) \ (\text{for more details, refer to \cite{M}}).  
\end{align*}
Let $Z(\mathfrak{g}) = U(\mathfrak{g})^{\mathfrak{g}} = \{x \in U(\mathfrak{g}) \ | \ \text{ad}y.x=0 \ \forall \ y \in U(\mathfrak{g}) \}$ where $Z(\mathfrak{g})$ denotes the center of $U(\mathfrak{g})$. It is well-known that $Z(\mathfrak{g})$ contains only even elements (see (2.6) of \cite{CW}).

(4.5) We now define Gelfand invariants which are elements of $Z(\mathfrak{g})$.\\
Fix $k \in \mathbb{N}$ and define
\begin{align*}
\widetilde{T}_k = \sum_{(i_1, \cdots, i_k)} (-1)^{p(i_2) + \cdots + p(i_k)}E_{i_1 i_2} \otimes E_{i_2 i_3} \otimes \cdots \otimes E_{i_k i_1} 
\end{align*}  
where the summation runs over all $(i_1, \cdots, i_k) \ (1 \leqslant i_j \leqslant M+N)$. Certainly $\widetilde{T}_k \in T^k(\mathfrak{g}) = \mathfrak{g} \otimes \cdots \otimes \mathfrak{g} \ (k$ times). Now $\mathfrak{g}$ acts naturally on $T^k(\mathfrak{g})$ and we shall denote this action by ad. It can be directly verified that ad$x.\widetilde{T}_k=0$.\\
Let $T_k = \pi(\widetilde{T}_k)$ and $T_k \in Z(\mathfrak{g})$. The following result is well-known.

\bppsn
Let $T$ be the algebra generated by $T_k$ where $k \in \mathbb{Z}_{+}$. Then we have $T = Z(\mathfrak{g})$. 
\eppsn      

(4.6) We now extend the notions to $\mathfrak{g}
\otimes B$ where $B$ is an associative, commutative and finitely generated algebra with unity.\\
Let $T(\mathfrak{g} \otimes B) = \bigoplus_{k \in \mathbb{Z}_{+}}T^k(\mathfrak{g} \otimes B)$ be the corresponding tensor algebra for $\mathfrak{g} \otimes B$ and $U(\mathfrak{g} \otimes B)$ be its universal enveloping algebra. Then we have a natural surjective algebra homomorphism
\begin{align*}
\pi(B) : T(\mathfrak{g} \otimes B) \longrightarrow U(\mathfrak{g} \otimes B).    
\end{align*}  
Fix $a_1, \cdots, a_k \in B$. For $k \in \mathbb{N}$, let us define $\widetilde{T}_k(a_1, \cdots, a_k)$  
\begin{align*}
= \sum_{(i_1, \cdots, i_k)} (-1)^{p(i_2) + \cdots + p(i_k)}E_{i_1 i_2}(a_1) \otimes E_{i_2 i_3}(a_2) \otimes \cdots \otimes E_{i_k i_1}(a_k). 
\end{align*}  
As we have seen earlier $T^k(\mathfrak{g} \otimes B)$ is a $(\mathfrak{g} \otimes B)$-module and in particular a $\mathfrak{g}$-module. We shall denote this action by ad.
\blmma \label{L4.2}
ad$x.\widetilde{T}_k(a_1, \cdots, a_k)=0 \ \forall \ x \in \mathfrak{g}$ and $a_1, a_2, \cdots, a_k \in B, \ k \in \mathbb{N}$.    
\elmma
\begin{proof}
Fix $a_1, \cdots, a_k \in B$ and consider the vector space $\mathfrak{g} \otimes \mathbb{C}a_i = \mathfrak{g}(a_i)$.\\
Then $\mathfrak{g}(a_1) \otimes \cdots \otimes \mathfrak{g}(a_k)$ is a $\mathfrak{g}$-module and isomorphic to $\mathfrak{g} \otimes \cdots \otimes \mathfrak{g}$. Since $\widetilde{T}_k(a_1, \cdots, a_k)$ goes to $\widetilde{T}_k$, the lemma follows by (4.5).   
\end{proof}    
(4.7) Let us now put $T_k(a_1, \cdots, a_k) = \pi(B)\widetilde{T}_k(a_1, \cdots, a_k)$ 
along with\\ $U(\mathfrak{g} \otimes B)^{\mathfrak{g}} = \{x \in U(\mathfrak{g} \otimes B) \ | \ \text{ad}y.x=0 \ \forall \ y \in U(\mathfrak{g}) \}$.\\
Further note that $T_k(a_1, \cdots, a_k) \in U(\mathfrak{g} \otimes B)^{\mathfrak{g}}$.\\
As $\mathfrak{g}$ acts on $\mathfrak{g} \otimes B$ only on the first factor, we suggest the following.
\begin{conj}
Let $T(B)$ be the algebra generated by $T^k(a_1, \cdots, a_k)$ where $a_1, a_2, \cdots, a_k \in B$ and $k \in \mathbb{Z}_{+}$. Then $T(B) = U(\mathfrak{g} \otimes B)^{\mathfrak{g}}$.   
\end{conj}

(4.8) In the rest of this section, we consider $B = \mathbb{C}[t,t^{-1}]$ which is the Laurent polynomial algebra in one variable. For any given $\mathbb{Z}_2$-graded vector space $V = V_{\overline{0}} \oplus V_{\overline{1}}$, let us denote $L(V) = V \otimes B$. $L(V)$ is a $\mathbb{Z}_2$-graded vector space with even part $L(V_{\overline{0}})$ and odd part $L(V_{\overline{1}})$. For $v \in V$ and $n \in \mathbb{Z}$, set $v(n)=v \otimes t^n$. We now recall the notion of evaluation modules for $L(\mathfrak{g})$. We shall closely follow Section $3$ of \cite{R2}.

(4.9) \textbf{Evaluation modules for $L(\mathfrak{g})$.} Fix $n \in \mathbb{N}$ and let $V(\lambda_1), \cdots, V(\lambda_n)$ be finite dimensional irreducible highest weight modules for $\mathfrak{g}$ with highest weights $\lambda_1, \cdots, \lambda_n$. Suppose that $d_1, \cdots, d_n$ are non-zero distinct complex numbers and put $\underline{d}=(d_1, \cdots, d_n) \in \mathbb{C}^n$. Let us now consider 
\begin{align*}
V(\underline{\lambda}, \underline{d}) = \bigotimes_{i=1}^{n} V(\lambda_i).
\end{align*}
Then $V(\underline{\lambda}, \underline{d})$ can be made into a $L(\mathfrak{g})$-module in the following way.\\
Let $x \in \mathfrak{g}$ and $w_i \in V(\lambda_i)$ be homogeneous elements. Define
\begin{align*}
x \otimes t^m.(w_1 \otimes \cdots \otimes w_n)\\
= \sum_{k=1}^{n} {d_k}^{m} (-1)^{|x|(|w_1| + \cdots + |w_{k-1}|)}w_1 \otimes \cdots \otimes x.w_k \otimes \cdots \otimes w_n.          
\end{align*}
It is an easy exercise to check that $V(\underline{\lambda}, \underline{d})$ is a $L(\mathfrak{g})$-module. In fact, we shall now see that it is also irreducible. Consider the evaluation map
\begin{align*}
\pi(\underline{d}) : L(\mathfrak{g}) \longrightarrow \bigoplus \mathfrak{g} \ 
(n \ \text{copies})
\end{align*}
$\,\,\,\,\,\,\,\,\,\,\,\,\,\,\,\,\,\,\,\,\,\,\,\,\,\,\,\,\,\,\,\,\,\,\,\,\,\,\,\,\,\,\,\,\,\,\,\,\,\,\,\,\,\,\,\,\,\,\,\,\,\,\,\,\,\,\,\ x \otimes t^m \longmapsto ({d_1}^{m}x, {d_2}^{m}x, \cdots, {d_n}^{m}x)$.\\
Since $d_1, \cdots, d_n$ are non-zero distinct complex numbers, it can be shown that $\pi(\underline{d})$ is surjective by means of a Vandermonde argument. Now $V(\underline{\lambda}, \underline{d})$ is clearly an irreducible module over $\bigoplus \mathfrak{g}$ and the $L(\mathfrak{g})$ action is nothing but the lift from $\pi(\underline{d})$. This proves that
$V(\underline{\lambda}, \underline{d})$ is an irreducible $L(\mathfrak{g})$-module.\\
For each $1 \leqslant i \leqslant n$, define $p_i(t) \in \mathbb{C}[t,t^{-1}]$ by setting    
\begin{align*}
p_i(t) = \prod_{j \neq i} \dfrac{t-d_j}{d_i-d_j} \ \text{and} \ p(t) = \prod_{j=1}^{n} (t-d_j).
\end{align*}
Then it can be trivially verified that $p_i(d_j)= \delta_{ij}$ and $\sum_{i=1}^{n} p_i(t)=1$.\\
Let $J$ be the ideal generated by $p(t)$ inside $\mathbb{C}[t,t^{-1}]$. But again since we have $p(d_i)=0 \ \forall \ 1 \leqslant i \leqslant n$, it directly follows that $(\mathfrak{g} \otimes J). V(\underline{\lambda}, \underline{d})= (0)$ and Ker$(\pi(\underline{d}))= \mathfrak{g} \otimes J$. It is also easy to see that $\{p_i(t) \}_{i=1}^{n}$ forms a basis of the vector space $\mathbb{C}[t,t^{-1}]/J$.

(4.10) Consider the operator $T_k(a_1, \cdots, a_k)$ which acts trivially on $V(\underline{\lambda}, \underline{d})$ if any one of the $a_i$'s belongs to $J$. We also note that $T_k(a_1, \cdots, a_k)$ is linear in all the variables and therefore the only operators that are of interest are    
\begin{align*}
T_k \big (p_{j_1}(t), \cdots, p_{j_k}(t) \big ) \ \text{where} \ \{j_1, \cdots, j_k \} \subseteq \{1, \cdots, n \}.
\end{align*}
For simplicity, we shall denote the above operator by $T_k(j_1, \cdots, j_k)$. 

(4.11) Consider the weight space decomposition of $V(\underline{\lambda}, \underline{d})$ given by
\begin{align*}
V(\underline{\lambda}, \underline{d}) = \bigoplus_{\mu \in \mathfrak{h}^{*}} V(\underline{\lambda}, \underline{d})_{\mu}
\end{align*} 
where $V(\underline{\lambda}, \underline{d})_{\mu} = \{v \in V(\underline{\lambda}, \underline{d}) \ | \ h.v = \mu(h)v \ \forall \ h \in \mathfrak{h} \}$.\\
Denote the set of all $\mu$-highest weight vectors lying inside $V(\underline{\lambda}, \underline{d})$ by 
\begin{align*}
V(\underline{\lambda}, \underline{d})_{\mu}^{+} = \{v \in V(\underline{\lambda}, \underline{d})_{\mu} \ | \ \mathfrak{g}^{+}.v=(0) \}.
\end{align*} 
As $T_k(a_1, \cdots, a_k)$ commutes with the action of $\mathfrak{g}$, it follows that
\begin{align*}
T_k(j_1, \cdots, j_k).V(\underline{\lambda}, \underline{d})_{\mu}^{+} \subseteq V(\underline{\lambda}, \underline{d})_{\mu}^{+}. 
\end{align*}
This is a definitive way of producing more highest weight vectors and will be useful in decomposing the tensor product module $V(\underline{\lambda}, \underline{d})$. Several examples have been worked out for the Lie algebra case in \cite{R2}. In the setup of Lie algebras, the action of $T(\mathbb{C}[t,t^{-1}])$ on $V(\underline{\lambda}, \underline{d})_{\mu}^{+}$ is irreducible (see \cite{Ku,R2}). 
\brmk
$T_k = T_k(1, \cdots, 1)$.
\ermk
The above $T_k$ is the classical Gelfand invariant and acts as a scalar on each isotypical component of $V(\underline{\lambda}, \underline{d})$. Now since we already know that $\sum_{i=1}^{n}p_i(t)=1$, it is evident that
\begin{align*}
T_k = \sum_{(j_1, \cdots, j_k)} T_k(j_1, \cdots, j_k). 
\end{align*}
Although each $T_k$ acts as a scalar, it decomposes into several operators which do not act by scalars on an isotypical component.

(4.12) We shall end this section with yet another interesting operator which will be utilized in decomposing $\mathfrak{g}$-modules with their even parts.\\
Recall that $\mathfrak{g}_{\overline{0}}=($End$V)_{\overline{0}}$ which is a reductive Lie algebra.\\ 
For $k \in \mathbb{N}$, define
\begin{align*}
S_k = \sum_{(i_1, \cdots, i_k)}E_{i_1 i_2}E_{i_2 i_3} \cdots E_{i_k i_1} 
\end{align*}  
where $(i_1, \cdots, i_k)$ runs over all possible indices satisfying $1 \leqslant i_j \leqslant M+N$.

\bppsn
$[x,S_k]=0 \ \forall \ x \in \mathfrak{g}_{\overline{0}}$.
\eppsn
\begin{proof}
Since we are taking superbracket with respect to the even part, the required calculation takes place only in the setting of general linear Lie algebra where this calculation is standard (see \cite{R2}).
\end{proof}
Now let $V$ be a finite dimensional irreducible module for $\mathfrak{g}$. Since $\mathfrak{g}_{\overline{0}}$ is reductive and the additional center is contained in $\mathfrak{h}$, this center acts via semisimple operators on $V$. Consequently $V$ is a completely reducible module over $\mathfrak{g}_{\overline{0}}$. The operators $S_k$ will be essential in finding $\mathfrak{g}_{\overline{0}}$ highest weight vectors in $V$.

\section{Construction Of Anti-invariants}\label{S5}
In this section we work with general Lie superalgebras. The notation are in continuation of the previous section.

(5.1) Fix $\mathfrak{g} = \mathfrak{gl}(M,N)$ where $M+N \geqslant 3$. Take I $= \{1,2, \cdots, M+N \}$ along with I$_0= \{1,2, \cdots, M \}$ and I$_1= \{M+1, M+2, \cdots, M+N \}$.\\
We now define parity for I by setting $p(i) = \overline{0}$ if $i \in $ I$_0$ and  $p(i) = \overline{1}$ if $i \in $ I$_1$ so that $|\alpha_{ij}| = p(i) + p(j)$.\\
Recall that $Z(\mathfrak{g}) = U(\mathfrak{g})^{\mathfrak{g}} = \{x \in U(\mathfrak{g}) \ | \ \text{ad}y.x=0 \ \forall \ y \in \mathfrak{g} \}$.\\
In the previous section, we have defined $T_k \in U(\mathfrak{g})$ and noted that they actually generate $Z(\mathfrak{g})$ as an algebra.

(5.2) In \cite{G}, Gorelik defined anticenter of $U(\mathfrak{g})$ which we shall now recall. Define twisted ad-action of $U(\mathfrak{g})$ by setting 
\begin{align*}
\text{ad}^{\prime}g.u = gu - (-1)^{|g|(|u| + 1)}ug
\end{align*}
for all homogeneous elements $g \in \mathfrak{g}$ and $u \in U(\mathfrak{g})$.

\bdfn 
$\mathcal{A}(\mathfrak{g}) = \{u \in U(\mathfrak{g}) \ | \ \text{ad}^{\prime}g.u=0 \ \forall \ g \in \mathfrak{g} \}$. 
\edfn
We call $\mathcal{A}(\mathfrak{g})$ the anticenter of $U(\mathfrak{g})$ and refer to its elements as anti-invariants. 

\blmma \cite[(2.1)]{G}\label{L5.2}
For homogeneous elements $g \in \mathfrak{g}, \ m_1, m_2 \in U(\mathfrak{g}),\\  \mathrm{ad}^{\prime}g.(m_1m_2) = (\mathrm{ad}^{\prime}g.m_1)m_2 + (-1)^{|g|(|m_1| + 1)} m_1(\mathrm{ad}g.m_2)$.
\elmma      

In this section, we define $D_l \ (l \in \mathbb{N})$ and verify that $D_l \in \mathcal{A}(\mathfrak{g})$.\\
Define $D_l = \sum_{(i_1, \cdots, i_l)} (-1)^{p(i_1, i_2 \cdots, i_l)}E_{i_1i_2}E_{i_2i_3} \cdots E_{i_li_1}$\\
where the above summation runs over all possible $l$-tuples $(i_1, i_2, \cdots, i_l)$ satisfying $1 \leqslant i_k \leqslant M+N$. Here $p(i_1, i_2 \cdots, i_l) = p(i_1) + p(i_2) + \cdots + p(i_l)$.\\
We also define $p_r(i_1, i_2 \cdots, i_l) = p(i_1) + \cdots + \widehat{p}(i_r) + \cdots + p(i_l)$.

\bppsn
$D_l \in \mathcal{A}(\mathfrak{g}) \ \forall \ l \in \mathbb{N}$.
\eppsn    
\begin{proof}
Consider $E_{j_1j_2} \in \mathfrak{g}$. Then by Lemma \ref{L5.2}, we have 
\begin{align*}
\mathrm{ad}^{\prime}E_{j_1j_2}.D_l = \sum_{(i_1, \cdots, i_l)} (-1)^{p(i_1, i_2 \cdots, i_l)} \big ( (E_{j_1j_2}.E_{i_1i_2}) - \\ (-1)^{|\alpha_{j_1j_2}|(|\alpha_{i_1i_2}| + 1)} (E_{i_1i_2}.E_{j_1j_2}) \big ). 
E_{i_2i_3} \cdots E_{i_li_1} + \\  \sum_{(i_1, \cdots, i_l)} (-1)^{p(i_1, i_2 \cdots, i_l) + |\alpha_{j_1j_2}|(|\alpha_{i_1i_2}| + 1)}E_{i_1i_2} (\mathrm{ad}E_{j_1j_2}. E_{i_2i_3} \cdots E_{i_li_1})  
\end{align*}
(note that $(E_{j_1j_2}.E_{i_1i_2})$ and $(E_{i_1i_2}.E_{j_1j_2})$ are matrix multiplications)
  
(5.3) \,\,\,\,\,\,\,\,\,\,\,\,\,\,\,\,\,\,\,\ $= \sum_{r=1}^{l} (B_r - B_{r}^{\prime})$ where for each $2 \leqslant r \leqslant l$, we set \\ $B_r =  \sum_{(i_1, \cdots, i_l), i_r=j_2} (-1)^{p(i_1, i_2, \cdots, i_l) + |\alpha_{j_1j_2}|(|\alpha_{i_1i_2}| + 1)} \\ \epsilon(\alpha_{j_1j_2}, \sum_{k=2}^{r-1}\alpha_{i_ki_{k+1}}) E_{i_1i_2} \cdots E_{i_{r-1}j_2}E_{j_1i_{r+1}} \cdots E_{i_li_1}$ and\\
$B_{r}^{\prime} =  \sum_{(i_1, \cdots, i_l), i_{r+1}=j_1} (-1)^{p(i_1, i_2 \cdots, i_l) + |\alpha_{j_1j_2}|(|\alpha_{i_1i_2}| + 1)} \\ \epsilon(\alpha_{j_1j_2}, \sum_{k=2}^{r}\alpha_{i_ki_{k+1}}) E_{i_1i_2} \cdots E_{i_rj_2}E_{j_1i_{r+2}} \cdots E_{i_li_1}$\\
$(\text{here} \ i_{l+1} \ \text{is understood} \ as \ i_1 \ \text{and} \ \epsilon(\alpha, \beta) = (-1)^{|\alpha||\beta|})$.\\   
Also let $B_1 =  \sum_{(i_1, \cdots, i_l), i_1=j_2} (-1)^{p(i_1, \cdots, i_l)} E_{j_1i_2}E_{i_2i_3} \cdots E_{i_lj_2}$ and\\
$B_{1}^{\prime} =  \sum_{(i_1, \cdots, i_l), i_{2}=j_1} (-1)^{p(i_1, i_2 \cdots, i_l) + |\alpha_{j_1j_2}|(|\alpha_{i_1i_2}| + 1)}$.\\
$Claim \ 1.$ $B_r = B_{r-1}^{\prime} \ \forall \ 3 \leqslant r \leqslant l$,\\
$Claim \ 2.$ $B_2 = B_1^{\prime}$ and
$Claim \ 3.$ $B_1 = B_l^{\prime}$.\\
Once we prove these claims, the proposition will follow from (5.3).\\
$Proof \ of \ Claim \ 1.$ We can easily see that the main terms are the same for both $B_r$ and $B_{r-1}^{\prime}$. We only need to show that the signs are the same.\\
Fix $(i_1, \cdots, i_l)$ such that for $B_r, \ i_r=j_2$ and for $B_{r-1}^{\prime}, \ i_r=j_1$.\\
All other indices are the same. Clearly there is an one-to-one correspondence between those indices.\\
Now the parity for $B_r$ and $B_{r-1}^{\prime}$ are respectively given by\\  
$p_r(i_1, \cdots, i_l) + p(j_2) + |\alpha_{j_1j_2}|(|\alpha_{i_1i_2}| + 1) + |\alpha_{j_1j_2}||\sum_{k=2}^{r-1}\alpha_{i_ki_{k+1}}|$ and\\         
$p_r(i_1, \cdots, i_l) + p(j_1) + |\alpha_{j_1j_2}|(|\alpha_{i_1i_2}| + 1) + |\alpha_{j_1j_2}||\sum_{k=2}^{r-1}\alpha_{i_ki_{k+1}}|$.\\
We need to show that the parity is same.\\
There are several terms which are common and after cancelling them we need to prove that
$p(j_2) + |\alpha_{j_1j_2}||\alpha_{i_{r-1}j_2}| = p(j_1) + |\alpha_{j_1j_2}||\alpha_{i_{r-1}j_1}|$.\\
This is easy to verify.\\ 
$Proof \ of \ Claim \ 2.$ Again we need to show that the signs are the same.\\
The parity of $B_2$ and $B_1$ are respectively given by\\
$p_2(i_1, i_2, \cdots, i_l) + p(j_2) + |\alpha_{j_1j_2}|(|\alpha_{i_1j_2}| + 1)$ and\\
$p_2(i_1, i_2, \cdots, i_l) + p(j_1) + |\alpha_{j_1j_2}|(|\alpha_{i_1j_1}| + 1)$.\\
After cancelling the common terms, we can check that the signs match.\\
$Proof \ of \ Claim \ 3.$ Here again we need to verify the signs.\\
The parity of $B_1$ is $p_1(i_1, i_2, \cdots, i_l) + p(j_2)$ and the parity of $B_l$ is\\ $p_1(i_1, i_2, \cdots, i_l) + p(j_1) + |\alpha_{j_1j_2}|(|\alpha_{j_1i_2}| + 1) + |\alpha_{j_1j_2}||\sum_{k=2}^{l}\alpha_{i_ki_{k+1}}|$.\\
Note that $|\alpha_{j_1j_2}||\sum_{k=2}^{l}\alpha_{i_ki_{k+1}}| = |\alpha_{j_1j_2}||\alpha_{i_1i_2}|$ as $\sum_{k=1}^{l}\alpha_{i_ki_{k+1}}$ is in fact zero and in particular even. Also observe that $i_1=j_1$.\\
After cancelling the common terms, we need to prove that\\
$p(j_2) = p(j_1) + |\alpha_{j_1j_2}| = p(j_1) + p(j_1) + p(j_2)$ which is clearly true.\\
This completes the proofs of all the claims. 

\brmk
We shall record the following special case.
\begin{enumerate}
\item $\sum_{(i_1, i_2)} (-1)^{p(i_1) + p(i_2)} E_{i_1i_2}E_{i_2i_1}$ is an anti-invariant.
\item $\sum_{(i_1, i_2)} (-1)^{p(i_2)} E_{i_1i_2}E_{i_2i_1}$ is an invariant. It can be easily shown that it is infact a Casimir operator. Here we take the non-degenerate invariant bilinear form as supertrace form. We have three different expressions for the Casimir operator in case of $\mathfrak{g} = \mathfrak{gl}(M,N)$. The first two are defined in (2.11) and Proposition \ref{P2.3}. They have been proved to be equal in (2.15).          
\end{enumerate}
\ermk               
\end{proof}

\end{document}